\documentclass[12pt]{article}

\usepackage{amsmath}
\usepackage{amssymb}

\addtocounter{section}{-0}
\newtheorem{guess}{Theorem}[section]

\newtheorem{lemma}[guess]{Lemma}

\numberwithin{equation}{section}

\begin{document}

\title{On nonlinear interpolation}

\author{T. Kappeler\footnote{Supported in part by
the Swiss National Science Foundation.}, P. Topalov\footnote{Supported
in part by the NSF grant DMS-0901443.}}

\maketitle

\date{}

\begin{abstract}
\noindent In a case study on asymptotics of spectral quantities of Schr\"odinger operators we show how 
the Riesz-Thorin theorem on the interpolation of linear operators can be extended to nonlinear maps. 

\vspace{6pt}

\end{abstract}


\section{Introduction}
\label{1. Introduction}

Instead of stating a general abstract theorem on nonlinear interpolation we prefer
to focus on how such a result can be applied in a fairly simple context. The
arguments of the proofs can easily be adapted to different set-ups.

\medskip

For any $0 \leq s < \infty$ denote by $H^s_{\mathbb C} \equiv H^s({\mathbb T},
{\mathbb C})$ the complex Sobolev space of order $s$ of the $1d$ torus 
${\mathbb T} = {\mathbb R} / {\mathbb Z}$. Its elements can be represented by a Fourier
series $\varphi (x) = \sum _{k \in {\mathbb Z}} \hat \varphi (k) e^{2\pi i k x}$, where
$\hat \varphi (k) := \int ^1_0 \varphi (x) e^{-2\pi i k x} dx \in {\mathbb C}$, and
   \[ \| \varphi \| _s := \left( \sum _{k \in {\mathbb Z}} (1 + |k|)^{2s} |\hat
      \varphi (k)|^2 \right) ^{1/2} < \infty .
   \]
Note that $H^{s'}_{\mathbb C} \hookrightarrow H^s_{\mathbb C}$ for any $0 \le s \leq s'$. 
Denote by $B^s_{\mathbb C}(R) \subseteq H^s_{\mathbb C}$ the open ball of radius $R$, 
centered at $0$. Furthermore, for any $1 \leq p < \infty$ and $t \in {\mathbb R}$
denote by $\ell^{p,t}_{\mathbb C} \equiv \ell ^{p,t}(\mathbb{Z}, \mathbb C)$ the complex sequence space
   \[ \ell ^{p,t}_{\mathbb C} = \left\{ \xi = (\xi _j)_{j \in\mathbb{Z}} \subseteq
      {\mathbb C} \big\arrowvert \| \xi \| _{p,t} < \infty \right\}
   \]
where
   \[ \| \xi \| _{p,t} := \left( \sum _{j \in\mathbb{Z}} \big((1 + |j|)^t |\xi _j|\big)^p \right)^{1/p} .
   \]
Note that $\ell ^{p,t'}_{\mathbb C} \hookrightarrow \ell ^{p,t}_{\mathbb C}$ for
any $0 \leq t \leq t'$. Finally, for $b>a\ge 0$ and $\alpha\ge 0$, $\beta>0$, 
assume that $F : B^a_{\mathbb C}(R) \rightarrow \ell ^{p,\alpha + \beta a}_{\mathbb C}$ is an analytic map, 
bounded by $M_a > 0$, so that for some $M_b > 0$,
   \begin{equation}
    \label{1.1}     \sup _{\varphi \in B^b_{\mathbb C}(R)} \| F(\varphi )\|
                    _{p,\alpha + \beta b} \leq M_b\,.
    \end{equation}
By the characterization of analytic maps with sequence spaces as their range
(cf. e.g. \cite[Appendix A]{KP}) and the analyticity of $F$ it follows from
\eqref{1.1} that $F| _{B^b_{\mathbb C}(R)} : 
B^b_{\mathbb C}(R)\rightarrow \ell ^{p,\alpha + \beta b}_{\mathbb C}$
is analytic as well. Hence, equivalently, we are given the following commutative diagram
\begin{equation}\label{eq:diagram}
\begin{array}{ccc}
B^a_{\mathbb C}(R)&\stackrel{F}{\longrightarrow}&\lefteqn{\ell^{p,\alpha + \beta a}_{\mathbb C}}\\
\uparrow&  &\uparrow\\
B^b_{\mathbb C}(R)&\stackrel{F| _{B^b_{\mathbb C}(R)}}{\longrightarrow}&
\lefteqn{\ell^{p,\alpha + \beta b}_{\mathbb C}}
\end{array}
\end{equation}
where the horizontal maps $F$ and $F| _{B^b_{\mathbb C}(R)}$ are analytic and bounded by
the constants $M_a>0$ and $M_b>0$ respectively, and the vertical arrows denote the standard inclusions
of the corresponding spaces.

\medskip

The following theorem is an instance of an extension of the Riesz-Thorin theorem
(cf. e.g. \cite{BL}) to nonlinear maps and is inspired by a special case treated
in \cite[Appendix 2]{Ku}.

\medskip

\begin{guess}\label{th:interpolation} 
Under the assumptions stated above, $F$ interpolates between
$B^a_{\mathbb C}(R)$ and $B^b_{\mathbb C}(R)$. More precisely, for any
$s= (1 - \lambda )a + \lambda b$ with $0 \leq \lambda \leq 1$,
$F(B^s_{\mathbb C}(R)) \subseteq \ell ^{p,\alpha + \beta s}_{\mathbb C}$ and 
$F|_{B^s_{\mathbb C}(R)} : B^s_{\mathbb C}(R)\rightarrow \ell^{p,\alpha + \beta s}_{\mathbb C}$
is an analytic map, bounded by  $(M_a)^{1 - \lambda } (M_b)^{\lambda }$.
\end{guess}

\medskip

The proof of Theorem~\ref{th:interpolation}, presented in Section 2, uses in a crucial
way the assumption that $F$ is analytic on a ball in a complex (Hilbert) space.
The proof does not apply to nonlinear maps
   \[ F : B^s_{\mathbb C}(R) \cap H^s({\mathbb T}, {\mathbb R}) \rightarrow
      \ell ^{p, \alpha + \beta s} (\mathbb{Z}, {\mathbb R})
   \]
which are merely real analytic.
In the sequel we would like to discuss a type of problems where nevertheless
Theorem~\ref{th:interpolation} can be applied to real analytic nonlinear maps.
For $q$ in $L^2_0({\mathbb T}, {\mathbb R}) = \big\{ q \in L^2({\mathbb T},
{\mathbb R}) \big\arrowvert \int ^1_0 q(x)dx = 0\big\}$ let $L(q)$ be the
Schr\"odinger operator $L(q):= - d^2_x + q$. The Dirichlet spectrum of
$L(q)$, considered on the interval $[0,1]$, is real and consists of
simple eigenvalues. We list them in increasing order $\mu _1 < \mu _2 <
\ldots $. Furthermore denote by $M(x, \lambda )$ the fundamental
solution of $L(q)$, i.e. the $2 \times 2$ matrix valued function
satisfying $L(q)M = \lambda M, \lambda \in {\mathbb C}$, and $M(0,\lambda )
= Id_{2 \times 2}$,
   \[ M(x,\lambda ) = \begin{pmatrix} y_1(x,\lambda ) &y_2(x, \lambda ) \\
      y'_1(x,\lambda ) &y'_2(x,\lambda ) \end{pmatrix}.
   \]
When evaluated at $\lambda = \mu _n$, the Floquet matrix $M(1,\lambda )$ is
lower triangular, hence its eigenvalues are given by $y_1(1, \mu _n)$ and
$y'_2(1, \mu _n)$. By the Wronskian identity, they satisfy $y_1(1, \mu _n)
y'_2(1, \mu _n) = 1$. By deforming the potential $q$ to the zero potential
along the straight line $tq, 0 \leq t \leq 1$, one sees that $(-1)^n y_1
(1, \mu _n) > 0$. Hence the Floquet exponents are given by $\pm {\kappa}
_n$ where
   \[ {\kappa}_n:= - \log \left( (-1)^n y_1(1,\mu _n) \right)
   \]
and $\log $ denotes the principal branch of the logarithm. It turns out that
together the $\mu _n$'s and ${\kappa}_n$'s form a system of canonical
coordinates for $L^2_0$ -- see \cite{PT}. The $\kappa '_n$'s play also an important
r\^ole for proving the property of $1$-smoothing of the periodic KdV equation
\cite{KST2}. We want to determine the asymptotics
of the ${\kappa}_n$'s as $n \rightarrow \infty $. To state them, introduce
for any $s \in {\mathbb R}_{\geq 0}$ the Sobolev spaces
   \[ H^s_0 \equiv H^s_0({\mathbb T}, {\mathbb R}) = H^s({\mathbb T},
      {\mathbb C}) \cap L^2_0({\mathbb T}, {\mathbb R})
   \]

and denote by $\langle \cdot , \cdot \rangle $ the standard inner product in
$L^2_0({\mathbb T}, {\mathbb R})$,
   \[ \langle f, g \rangle = \int ^1_0 f(x) g(x) dx .
   \]
In \cite{KST1}, the following theorem was proved.

\medskip

\begin{guess}
\label{Theorem 1.2} {\rm \cite{KST1} } For any $N \in {\mathbb Z}_{\geq 0}$,
   \[ {\kappa}_n = \frac{1}{2\pi n} \left( \langle q, \sin 2 \pi n x \rangle
      + \frac{1}{n ^{N + 1}} \ell ^2_n \right)
   \]
uniformly on bounded subsets of $H^N_0$. Here $\ell ^2_n$ denotes the
$n'$th component of a sequence in $\ell ^2 \equiv \ell ^2 ({\mathbb N},
{\mathbb R})$.
\end{guess}

\medskip

Our aim is to extend Theorem~\ref{Theorem 1.2} to any fractional order Sobolev
space $H^s_0$ with $s \in {\mathbb R}_{\geq 0}$.
To this end we need to extend various quantities to the complex Hilbert space
$H^s_{0,{\mathbb C}} := \big\{ q \in H^s_{\mathbb C} \big\arrowvert \int ^1_0q(x)\,dx = 0\big\} $. 
Let $L^2_{0,{\mathbb C}} \equiv H^0_{0,{\mathbb C}}$ and
denote by $B^s_{0,{\mathbb C}}(R)$ the complex ball 
$B^s_{0,{\mathbb C}}(R):= \big\{ q \in H^s_{0,{\mathbb C}} \big\arrowvert \| q \| < R \big\}$.
In section 3 we prove the following

\medskip

\begin{guess}
\label{Theorem 1.3} For any $R > 0$ there exits $n_R > 0$ so that for any
$n > n_R$, ${\kappa}_n$ can be extended analytically to $B^0_{0,{\mathbb C}}(R)$. 
Moreover, for any $N \in {\mathbb Z}_{\geq 0}$ the sequence
$({\kappa}_n)_{n > n_R}$ satisfies the estimate
   \[ {\kappa}_n = \frac{1}{2\pi n} \left( \langle q, \sin 2n \pi x
      \rangle + \frac{1}{n^{N + 1}} \ell ^2_n \right)
   \]
uniformly on $B^N_{0,{\mathbb C}}(R)$.
\end{guess}

\medskip

At the end of Section 3 we show that Theorem~\ref{Theorem 1.3} allows to apply 
Theorem~\ref{th:interpolation} to  generalize Theorem~\ref{Theorem 1.2} to any fractional order Sobolev
space. 


\medskip

\begin{guess}
\label{Theorem 1.4} For any $s \in {\mathbb R}_{\geq 0}$ 
   \[ \kappa _n = \frac{1}{2\pi n} \left( \langle q, \sin 2n\pi x \rangle +
      \frac{1}{n^{s+1}} \ell ^2_n \right)
   \]
uniformly on bounded subsets of $H^s_0$.
\end{guess}

\medskip

By the same method one can show that similar results hold for many other spectral
quantities. We state without proof another such result. For $q \in L^2_0$, the
periodic spectrum of $L(q)$ on the interval $[0,2]$ is real and discrete. When
listed in increasing order and with their multiplicities the eigenvalues satisfy
   \[ \lambda _0 < \lambda ^-_1 \leq \lambda ^+_1 < \lambda ^-_2 \leq
      \lambda ^+_2 < \ldots .
   \]
Using \cite[Theorem 1.3, Theorem 1.4]{KST1}, one then shows with the help 
of Theorem~\ref{th:interpolation} that the following result holds.

\medskip

\begin{guess}
\label{Theorem 1.5} For any $s \in {\mathbb R}_{\geq 0}$
   \[ \left( \lambda ^+_n + \lambda ^-_n \right) / 2 - \mu _n = \langle q,
      \cos 2\pi n x \rangle + \frac{1}{n^{s+1}} \ell ^2_n
   \]
uniformly on bounded subsets of $H^s_0$.
\end{guess}

\medskip

{\em Related work:} In the context of nonlinear PDEs, Tartar obtained a nonlinear interpolation
theorem in \cite{Ta} which later was slightly improved in \cite[Theorem 1]{BS}. 
However, for nonlinear maps such as the ones encountered in the
analysis of Schr\"odinger operators the assumptions of these theorems are not
satisfied and hence they cannot be applied.

\medskip

{\em Acknowledgment:} The motivation of this paper originated from our observation 
that Theorem 6.1 in the paper \cite{SaShk} of  A. Savchuk and A. Shkalikov, the  theorem of Tartar on nonlinear interpolation \cite{Ta} cannot be applied as stated by the authors. After having pointed this out to A. Shkalikov, he sent us the statement of an unpublished interpolation theorem obtained by them earlier,
which can be used to prove their results.  Due to unfortunate personal circumstances their paper with a proof of this interpolation theorem was delayed. In order to be able to further progress on our projects, we came up with our own version, tailored to our needs and with an optimal estimate on the bounds of the interpolated maps. After having sent our preprint to A. Shkalikov, he sent us a preprint with the proof of their theorem \cite{SaShk1}. As the two papers turned out to be quite different, we agreed to publish them independently.


\section{Proof of Theorem~\ref{th:interpolation}}
\label{Proof of th:interpolation}

The proof of Theorem~\ref{th:interpolation} uses arguments of the proof of the
Riesz-Thorin theorem on the interpolation of linear operators -- see e.g.
\cite{BL}. Key ingredient is Hadamard's three-lines theorem.

\medskip

{\em Proof of Theorem~\ref{th:interpolation}.} For any $k \in \mathbb{Z}$, denote by $F_k$
the $k$'th component of $F$. As $F : B^a_{\mathbb C}(R)\to\ell^{p,\alpha+\beta b}_{\mathbb C}$ is analytic, 
$F_k : B^a_{\mathbb C}(R)\rightarrow {\mathbb C}$ is analytic. In particular, the restriction 
$F_k|_{B^s_\mathbb{C}(R)}$ of $F_k$ to $B^s_{\mathbb C}(R)\subseteq B^a_{\mathbb C}(R)$ is analytic for
any $a \leq s\leq b$. Hence once the claimed bound for $F|_{B^s_{\mathbb C}(R)}$ is established, 
the analyticity of 
$F|_{B^s_{\mathbb C}(R)} : B^s_{\mathbb C}(R) \rightarrow \ell ^{p, \alpha + \beta s}_{\mathbb C}$ 
will follow (cf. e.g. \cite[Appendix A]{KP}). 
To prove the claimed bound consider $\varphi \in B^s_{\mathbb C}(R)$ and define the complex family
   \[ \varphi _z(x) := \sum _{k \in {\mathbb Z}} \langle k \rangle ^{s-z}
      \hat \varphi (k) e^{2\pi i k x}
   \]
where $z = u + iv \in {\mathbb C}$, $u,v\in\mathbb{R}$, and $\langle k \rangle := 1 + |k|$. 
As $u$ and $v$ are real we see that for any $\varphi\in B^s_\mathbb{C}(R)$,
   \begin{equation*}
      \| \varphi _{u + iv}\| ^2_u = \sum _{k \in {\mathbb Z}}
       \langle k \rangle ^{2u} \big\arrowvert \langle k \rangle
       ^{-u-iv+s} \hat \varphi (k) \big\arrowvert ^2 
       = \sum _{k \in {\mathbb Z}} \langle k \rangle ^{2s}
       |\hat \varphi (k)|^2 = \| \varphi \|^2_s < R^2 .
   \end{equation*}
In particular, for any $v\in\mathbb{R}$,
\begin{equation}\label{eq:F-ab-bounded}
\|F(\varphi _{a + iv})\| _{p,\alpha + \beta a}\le M_a\quad\mbox{and}\quad
\|F(\varphi _{b + iv})\| _{p,\alpha + \beta b}\le M_b\,.
\end{equation}
Similarly, for any $a \leq u \leq b$ and $v \in {\mathbb R}$,
   \[ \| \varphi _{u + iv}\| ^2_a = \sum _{k \in {\mathbb Z}}
       \langle k \rangle ^{2(a-u+s)} |\hat \varphi (k)|^2 
       \leq \sum _{k \in {\mathbb Z}} \langle k \rangle ^{2s}
       |\hat \varphi (k)|^2 =\|\varphi\|_s^2< R^2 .
   \]
Hence $\varphi _{u + iv} \in B^a_{\mathbb C}(R)$ for any $a \leq u \leq b$
and any $v \in {\mathbb R}$. In particular, $F(\varphi _{u+iv})\in\ell ^{p,\alpha +\beta a}_{\mathbb C}$ and
\begin{equation}\label{eq:F-bounded}
\|F(\varphi _{u+iv})\|_{p,\alpha+\beta a}\le M_a\,.
\end{equation}
As $\varphi _s = \varphi $ for any $\varphi \in B^s_{\mathbb C}(R)$, 
we need to prove the claimed estimate for $z = s$. 
Denote by $1 < q \leq \infty$ the number conjugate to $p$, i.e. $\frac{1}{p} + \frac{1}{q} = 1$, and
consider an arbitrary sequence $\xi = (\xi _k)_{k \in \mathbb{Z}} \subseteq{\mathbb C}$ with 
{\em finite support} so that $\| \xi \|_{q, - (\alpha + \beta s)} \leq 1$. 
Similarly as in the case of $H^s_{\mathbb C}$, define for $z = u + iv \in {\mathbb C}$, $u,v\in\mathbb{R}$,
   \[ \xi _z := \left( \langle k \rangle ^{-\beta (s - z)} \xi _k \right) 
      _{k \in \mathbb{Z}}. 
   \]
Then
   \begin{align}\label{eq:xi-ab-bounded}
                     \| \xi _{u+iv}\|^q_{q,-(\alpha + \beta u)} &=
                     \sum_{k \in \mathbb{Z}} \Big(\langle k\rangle^{-(\alpha + \beta u)}
                     \big|\langle k \rangle ^{-\beta(s - z)}\xi _k\big|\Big) ^q\nonumber\\
                  &\leq \sum _{k \in \mathbb{Z}} \left( \langle k \rangle ^{-(\alpha
                     + \beta s)} |\xi _k| \right) ^q \leq 1\,.
   \end{align}
Similarly as above, $\xi _s = \xi $, and for any $a \leq u \leq b$ and $v \in {\mathbb R}$,
\begin{equation}\label{eq:xi-bounded}
\|\xi_{u+iv}\|_{q,-(\alpha + \beta a)}\le\| \xi\|_{q,-\alpha + \beta b}\,.
\end{equation}
Denote by $\langle \cdot , \cdot\rangle $ the $\ell ^2$ dual pairing as well as its extension to $\ell ^{p,
\alpha + \beta a}_{\mathbb C} \times \ell ^{q,-(\alpha + \beta a)}_{\mathbb C}$
and introduce the vertical strip
   \[ S_{a,b} = \left\{ z = u + iv \big\arrowvert v \in {\mathbb R} , \ a <
      u < b \right\}
   \]
and its closure $\overline S_{a,b}$. To obtain the claimed estimate we want to
apply Hadamard's three lines theorem to the following function
   \begin{equation*}
             f : \overline S_{a,b} \rightarrow {\mathbb C}, \ z \mapsto
               \langle F(\varphi _z), \xi _z \rangle 
   \end{equation*}
where
  \[ \langle F(\varphi _z), \xi _z \rangle = \sum _{k \in \mathbb{Z}} F_k(\varphi _z)
      \langle k \rangle ^{-\beta (s-z)} \xi _k .
   \]
As the support of $\xi$ is finite, the latter sum is finite, and hence the function $f$ is well defined
and $\|\xi\|_{q,-\alpha+\beta b}<\infty$.
In view of \eqref{eq:F-bounded} and \eqref{eq:xi-bounded}, 
$f: \overline S_{a,b} \rightarrow {\mathbb C}$ is {\em bounded} as
for any $z\in \overline S_{a,b}$,
\[
|f(z)|=|\langle F(\varphi _z), \xi _z \rangle|\le
\|F(\varphi _z)\|_{p,\alpha+\beta a}\|\xi_z\|_{q,-(\alpha+\beta a)}\le M_a\|\xi\|_{q,-\alpha+\beta b}.
\]
Note that on the strip $S_{a,b}$, the curves 
$z \mapsto \xi _z \in \ell ^{q,-(\alpha + \beta a)}_{\mathbb C}$ and 
$z \mapsto \varphi _z \in B^a_{\mathbb C}(R) \subseteq H^a_{\mathbb C}$ are analytic. 
As $F_k : B^a_{\mathbb C}(R) \rightarrow {\mathbb C}$ is analytic for any $k \in \mathbb{Z}$ it
then follows that $f$, being a finite sum of analytic functions, is
analytic on $S_{a,b}$ and continuous on the closure $\overline S_{a,b}$.
Moreover, in view of \eqref{eq:F-ab-bounded} and \eqref{eq:xi-ab-bounded}, the following estimates hold for any
$v \in {\mathbb R}$,
   \[ |f(a + iv)| \leq \| F(\varphi _{a + iv})\| _{p,\alpha + \beta a}
      \| \xi _{a + iv}\| _{q,-(\alpha + \beta a)} \leq M_a
   \]
and
   \[ |f(b + iv)| \leq \| F(\varphi _{b + iv})\| _{p,\alpha + \beta b}
      \| \xi _{b + iv}\| _{q,-(\alpha + \beta b)} \leq M_b .
   \]
Hence we can apply Hadamard's three-lines theorem to $f$ (cf. e.g.
\cite[Appendix to IX,.4]{RS}) to conclude that, for $s = (1 - \lambda ) a +
\lambda b$ and any $v \in {\mathbb R}$
   \[ |f(s + iv)| \leq (M_a)^{1 - \lambda } (M_b)^{\lambda } .
   \]
In particular, for $z = s$ one has
   \[ |f(s)| = \big\arrowvert \langle F(\varphi ), \xi \rangle \big\arrowvert
      \leq (M_a)^{1 - \lambda } (M_b)^{\lambda }
   \]
where we took into account that $\varphi _s = \varphi$ and $\xi _s = \xi $.
As the sequences $\xi $ with finite support are dense in $\ell ^{q,-(\alpha + \beta s)} _{\mathbb C}$ 
the claimed estimate $\| F(\varphi )\| _{p,\alpha +\beta s}\leq(M_a)^{1-\lambda } (M_b)^{\lambda }$ 
follows from Hahn-Banach theorem and the fact that $\ell ^{q,-(\alpha + \beta s)} _{\mathbb C}$
is the dual space of $\ell ^{p,\alpha + \beta s} _{\mathbb C}$. 
\hspace*{\fill }$\square$

\vskip 1 cm 

\section{Proofs of Theorem~\ref{Theorem 1.3} and Theorem~\ref{Theorem 1.4}}
\label{3. Proofs of Theorem 1.3 und Theorem 1.4}

Before proving Theorem~\ref{Theorem 1.3} we need to make some preparatory
considerations. Note that for $q$ in $L^2_{0,{\mathbb C}}$ the operator $L(q)$
is no longer symmetric with respect to the $L^2$-inner product $\langle f,g
\rangle = \int ^1_0 f(x) \overline{g(x)}\,dx$. The Dirichlet spectrum is
still discrete, but the eigenvalues might be complex valued and multiple.
We list them according to their algebraic multiplicities and in lexicographic
ordering
   \[ \mu _1 \preccurlyeq \mu _2 \preccurlyeq \mu _3 \preccurlyeq \ldots
   \]
where for any complex numbers $a, b$
   \[ a \preccurlyeq b \mbox { iff } \left[ \mbox{Re } a < \mbox{Re } b \right]
      \mbox { or } \left[ \mbox{Re } a = \mbox{Re } b \ \mbox { and }
      \mbox{ Im } a \leq \mbox{ Im } b \right] .
   \]
By \cite{PT} there exists for any $R > 0$ an integer $m_R > 0$ so that
for any $q \in L^2_{0,{\mathbb C}}$ with $\| q\| := \langle q,q \rangle ^{1/2}
< R$ the Dirichlet eigenvalues $\mu _n \equiv \mu _n(q), n \geq 1$, satisfy
the following estimates
   \[ |\mu _n - n^2 \pi ^2| < \pi / 4 \quad \forall\,\,n > m_R
   \]
and
   \[ |\mu _n| < m^2_R \pi ^2 + \pi / 4 \quad \forall 1 \leq n \leq m_R .
   \]
In particular, for any $n > m_R$, $\mu _n$ is a simple Dirichlet eigenvalue
of $L(q)$ and hence analytic on the complex ball 
$B^0_{0,{\mathbb C}}(R)= \{ q \in L^2_{0,{\mathbb C}} \big\arrowvert \| q\| < R \} $. 
To see that $\kappa _n$ can be analytically extended to $B^0_{0,{\mathbb C}}(R)$ for
$n$ sufficiently large we first note that by \cite{PT}, the Floquet matrix
$M(1, \lambda , q)$ is analytic on ${\mathbb C} \times L^2_{0,{\mathbb C}}$ and, 
with $\lambda = \nu ^2$, $y_1(1, \nu ^2, q)$ satisfies the following estimate
   \begin{equation}
   \label{3.1} \big\arrowvert y_2(1, \nu ^2, q) - \cos \nu \big\arrowvert \leq
               \frac{1}{|\nu |} \exp \left( \big\arrowvert \mbox{Im } \nu
               \big\arrowvert + \| q \| \right) .
   \end{equation}
As for $n > m_R$, one has $|\mu _n - n^2 \pi ^2| < \pi / 4$, it follows that
$\nu _n \equiv \nu _n(q) = \sqrt[+]{\mu _n(q)}$ is well defined. Here
$\sqrt[+]{z}$ denotes the principal branch of the square root. Then for any
$n > m_R$ and $q \in B^0_{0,{\mathbb C}}(R)$,
 $\nu _n = n\pi \sqrt[+]{1 + z_n},$ where $\ z_n = \frac{\mu _n - n^2\pi ^2}{n^2 \pi ^2} $
satisfies $|z_n| \leq \frac{1}{4\pi n^2} \leq \frac{1}{10}$ and thus
 \[ \big\arrowvert \sqrt[+]{1 + z_n} - 1 \big\arrowvert = 
                     \Bigg\arrowvert \frac{(\sqrt[+]{1 + z_n}-1)(\sqrt[+]
                     {1 + z_n} + 1)}{\sqrt[+]{1 + z_n} + 1} \Bigg\arrowvert
                     \leq |z_n| 
\]
or
\[
                  |\nu _n - n\pi | \leq \frac{\pi }{4} \frac{1}{4n} \leq
                     \frac{1}{4n}; \quad |\mbox{Im } \nu _n | \leq \frac{1}{4n} .
  \]
These estimates are used in the asymptotics of $y_1(1, \nu ^2_n, q)$. As
$(-1)^n$ $\cos \nu _n = 1 + (\nu _n - n\pi ) \int ^1_0 - \sin (t(\nu _n - n\pi ))
dt$ one concludes that $|(-1)^n \cos \nu _n - 1| \leq \frac{1}{4n}$ for any
$n > m_R$ and as $|\nu _n| \geq n\pi - \frac{1}{4n} \geq 2n$ it follows from
\eqref{3.1} that
   \[ |(-1)^n y_1(1, \mu _n, q) - 1 | \leq \frac{1}{4n} + \frac{1}{2n} 
      \exp \left( \frac{1}{4n} + R \right)
   \]
for any $n > m_R$ and $\| q\| < R$. Now choose $n_R \geq m_R$ so large that
   \[ \frac{1}{4n} + \frac{1}{2n} \exp \left( \frac{1}{4n} + R \right) \leq
      \frac{1}{2} \quad \forall n > n_R .
   \]
As a consequence, for any $q \in L^2_{0,{\mathbb C}}$ with $\| q\| < R$ and
any $n > n_R$, 
   \[ \kappa _n = - \log \left( (-1)^n y_1(1, \mu _n, q) \right) 
   \]
is well defined. We thus have proved.

\medskip

\begin{lemma}
\label{Lemma 3.1} For any $R > 0$ there exists $n_R > 0$ so that $\forall
n > n_R$
   \[ \kappa _n :  B^0_{0,{\mathbb C}}(R) \rightarrow {\mathbb C} , \quad
      q \mapsto - \log \left( (-1) ^n y_1 (1, \mu _n, q) \right)
   \]
is well defined and analytic.
\end{lemma}

{\it Proof of Theorem~\ref{Theorem 1.3}.} In view of Lemma~\ref{Lemma 3.1}
it remains to prove for any $R > 0$, and any $N \in {\mathbb Z}_{\geq 0}$,
   \[ \sum _{n > n_R} n^{2N+2} \big\arrowvert 2\pi n \kappa _n - \langle
      q, \sin 2 n\pi x \rangle \big\arrowvert ^2
   \]
is bounded on $B^N_{0,{\mathbb C}}(R)$. Going through the arguments of the
proof of Theorem 1.1 in \cite{KST1}, one sees that this is indeed the case.
\hspace*{\fill }$\square $

\bigskip

It remains to show Theorem~\ref{Theorem 1.4}. As already mentioned in the
introduction, we will use the result of Theorem~\ref{th:interpolation} on 
nonlinear interpolation.

\medskip

{\it Proof of Theorem~\ref{Theorem 1.4}.} By Theorem~\ref{Theorem 1.2},
Theorem~\ref{Theorem 1.4} holds in the case where $s$ is an integer. Let
$N \in {\mathbb Z}_{\geq 0}$ and $N < s < N + 1$. For any given $R > 0$
choose $n_R > 0$ as in Theorem~\ref{Theorem 1.3}. It follows from
Theorem~\ref{Theorem 1.2} that there exists $M_R > 0$ such that for any
$q \in L^2_0$ with $\| q\| < R$ and any $1 \leq n \leq n_R$
   \[ \big\arrowvert 2\pi n \kappa _n - \langle q, \sin 2\pi n x \rangle
      \big\arrowvert \leq M_R
   \]
and therefore
   \begin{equation}
   \label{3.5} \left( \sum _{n \leq n_R} n^{2(s+1)} \big\arrowvert 2\pi n 
               \kappa _n - \langle q, \sin 2\pi n x \rangle \big\arrowvert
               ^2 \right) ^{1/2} \leq n^{s + 3/2}_R M_R .
   \end{equation}
On the other hand, by Theorem~\ref{Theorem 1.3}, for any $n > n_R, \kappa _n$
extends analytically to $B^0_{0,{\mathbb C}}(R)$ and there exist constants
$M_{N,R}, M_{N + 1,R}$ so that
   \[ \left( \sum _{n > n_R} n^{2(N+1)} \big\arrowvert 2\pi n \kappa _n -
      \langle q, \sin 2n\pi x \rangle \big\arrowvert ^2 \right) ^{1/2}
      \leq M_{N,R} \qquad \forall q \in B^N_{0,{\mathbb C}}(R)
   \]
as well as 
   \[ \left( \sum _{n > n_R} n^{2(N+2)} \big\arrowvert 2\pi n \kappa _n -
      \langle q, \sin 2n\pi x \rangle \big\arrowvert ^2 \right) ^{1/2}
      \leq M_{N + 1,R} \qquad \forall q \in B^{N +1}_{0,{\mathbb C}}(R).
   \]
Hence Theorem~\ref{th:interpolation} applies and implies that with $0 \leq
\lambda \leq 1$ chosen so that $s = (1 - \lambda )N + \lambda (N + 1)$,
   \[ \left( \sum _{n > n_R} n^{2(s + 1)} \big\arrowvert 2\pi n \kappa _n -
      \langle q, \sin 2n\pi x \rangle \big\arrowvert ^2 \right) ^{1/2}
      \leq M_{s,R} \quad \forall q \in B^s_{0,{\mathbb C}}(R)
   \]
where $M_{s,R} = (M_{N,R})^{(1-\lambda )}(M_{N+1,R})^\lambda $.
From this together with \eqref{3.5}, the claimed result follows.
\hspace*{\fill }$\square $


\end{document}